\documentclass[11pt, reqno]{amsart}
\usepackage{amsthm, amsmath, amsfonts, amssymb}
\usepackage{array, texdraw}
\evensidemargin \oddsidemargin\parindent=8mm

\begin{document}

\renewcommand{\proofname}{\bf Proof}
\newtheorem*{rem*}{Remark}
\newtheorem{cor}{Corollary}
\newtheorem{prop}{Proposition}
\newtheorem{lem}{Lemma}
\newtheorem{theo}{Theorem}
\newtheorem{fact}{Fact}
\newtheorem{hypo}{Hypothesis}
\newfont{\zapf}{pzcmi}

\def\R{{\mathbb R}}
\def\N{{\mathbb N}}
\def\Z{{\mathbb Z}}
\def\E{{\mathbb E}}
\def\P{{\mathbb P}}
\def\V{{\mathbb D}}
\def\I{{\mathbb I}}
\newcommand{\D}{\hbox{ \zapf D}}

\title[The area of exponential random walk]
{The area of exponential random walk \\and partial sums of uniform
order statistics}
\author{Vladislav V. Vysotsky}
\thanks{Partly supported by the grant NSh-4222.2006.1.}

\begin{abstract}
Let $S_i$ be a random walk with standard exponential increments.
We call $\sum_{i=1}^k S_i$ its {\it k-step area}. The random
variable $\inf \limits_{k \ge 1} \frac{2}{k(k+1)} \sum_{i=1}^k
S_i$ plays important role in the study of so-called
one-dimensional sticky particles model. We find the distribution
of this variable and prove that $$\P \biggl \{ \inf \limits_{k \ge
1} \frac{2}{k(k+1)} \sum_{i=1}^k S_i  \ge t \biggr \} = \P \biggl
\{ \inf \limits_{k \ge 1} \sum_{i=1}^k \bigl( S_i - it \bigr) \ge
0 \biggr \} = \sqrt{1-t} \, e^{-t/2}$$ for $0 \le t \le 1$. We
also show that
$$\lim \limits_{n \to \infty} \P \biggl \{ \min \limits_{1 \le k
\le n} \frac{2n}{k(k+1)} \sum_{i=1}^k U_{i, n} \ge t \biggr \} =
\sqrt{1-t} \, e^{-t/2} ,
$$ where $U_{i, n}$ are the order statistics of $n$ i.i.d. random
variables uniformly distributed on $[0,1]$.


% 60G50 Sums of independent random variables; random walks
% 62G30 Order statistics; empirical distribution functions
\end{abstract}
\maketitle

\section{Introduction}
Let $S_i$ be a positive random walk. We call $\sum_{i=1}^k S_i$,
where $k \ge 1$, its {\it k-step area}. We are interested in
$k$-step areas of an {\it exponential random walk} $S_i$, that is
a walk with standard exponential increments.

Let us normalize each $k$-step area dividing it by its expectation
$\E \sum_{i=1}^k S_i = \frac{k(k+1)}{2}$. Now introduce the random
variable $\inf \limits_{k \ge 1} \frac{2}{k(k+1)} \sum_{i=1}^k
S_i$, which is the main object of study in this paper. Our goal is
to find the probabilities
$$G(t) := \P \biggl \{ \inf \limits_{k \ge 1} \frac{2}{k(k+1)}
\sum_{i=1}^k S_i  \ge t \biggr \}.$$ We can also write the
right-hand side in a more pleasant form $$G(t) = \P \biggl \{ \inf
\limits_{k \ge 1} \sum_{i=1}^k \bigl( S_i - it \bigr) \ge 0 \biggr
\},$$ which resembles a ruin probability.

The function $G(t)$ arises in the study of so-called
one-dimensional {\it sticky particles model}. Let us briefly
describe the model. We consider a system of $n$ identical
particles, each one of mass $n^{-1}$. At time zero the immobile
particles are randomly distributed on the real line. The particles
begin to move under the forces of mutual attraction. When two or
more particles collide, they {\it stick} together forming a new
particle (``{\it cluster}'') whose characteristics are defined by
the laws of mass and momentum conservation. Between collisions
particles move according to the laws of Newtonian mechanics.

We suppose that the force of mutual attraction does not depend on
distance and equals the product of masses; this is very natural
for one-dimensional models. Thus at any moment, the acceleration
of a particle is equal to difference of masses to the right and to
the left of the particle.

There are two natural and well-known models of random initial
positions of particles. In the {\it uniform model}, $n$ particles
are uniformly and independently spread on $[0,1]$. In the {\it
Poisson model}, the particles are located at the points of first
$n$ jumps of a Poisson process with intensity $n$ (i.e., a
standard Poisson process multiplied by $n^{-1}$).

For more information about systems of sticky particles, see
\cite{Giraud}, \cite{LifShi}, \cite{MP}, \cite{I}, and references
therein.

Let us agree that the term ``particle'' refers to the initial
particles only, and let a ``cluster'' be a product of a collision
as well as an initial particle that has not experienced any
collisions. Thus at any moment $t>0$, the system consists of
clusters and each cluster contains one or more particles. As time
goes, clusters aggregate and became larger and larger while the
number of clusters decreases. Finally, at some moment all clusters
merge into a single cluster containing all (initial) particles.

By $K_n(t)$ denote the number of clusters at time $t$ in the
system of $n$ particles. This quantity is a random step function,
which decreases (in $t$) from its initial value $n$ to $1$. The
problem, which leads to the study of the function $G(t)$, is to
describe the asymptotics of $K_n(t)$ as $n \to \infty$. This
problem was introduced in \cite{I}, where the author proved the
following statement. Both in the uniform and the Poisson models of
initial positions, for any $t \ge 0$, we have
$$\frac{K_n(t)}{n} \stackrel{\D} {\longrightarrow}
K(t), \qquad n \to \infty,$$ where $K(t)$ is a deterministic
function satisfying $K(t) = e^{t^2} \bigl( G(t^2) \bigr)^2.$ It
was conjectured on the basis of numerical simulations that
$K(t)=1-t^2$ for $0 \le t \le 1$.

The main result of the current paper, the formula
\begin{equation}\label{main}
G(t)= \P \biggl \{ \inf \limits_{k \ge 1} \sum_{i=1}^k \bigl( S_i
- it \bigr) \ge 0 \biggr \} = \sqrt{1-t} \, e^{-t/2}, \qquad 0 \le
t \le 1,
\end{equation}
shows that the conjecture is true. Our study of the problem was
motivated  by the wish to prove the weird formula \eqref{main} as
well as by the necessity to verify some properties of $K(t)$ the
author needed in his further investigation of $K_n(t)$. New
results on the number of clusters will be soon published in
\cite{I2}.

We also note that in view of the well-know connection between
exponential random walks and order statistics, $G(t)$ could be
represented in the form
$$G(t) = \lim \limits_{n \to \infty} \P \biggl \{ \min \limits_{1
\le k \le n} \frac{2n}{k(k+1)} \sum_{i=1}^k U_{i, n}  \ge t \biggr
\}, $$ where $U_{i, n}$ are the order statistics of $n$ i.i.d.
random variables uniformly distributed on $[0,1]$. This equality
will be proved rigorously in the end of Sec.~\ref{Sec DE}. Thus
$G(t)$ is closely related to partial sums of uniform order
statistics.

The proof of \eqref{main} is organized as follows. In
Sec.~\ref{Sec Diff}~and~\ref{Sec Prop} we study properties of the
functions
$$G_n(t):= \P \biggl \{ \min \limits_{1 \le k \le n }
\frac{2}{k(k+1)} \sum_{i=1}^k S_i \ge t \biggr \},$$ which
converge to $G(t)$. We show that $G_n(t)$ are continuously
differentiable and that $G_n'(t)$ converge uniformly;
consequently, $G(t)$ has a continuous derivative. We also obtain
an ordinary differential equation for $G(t)$, but the right-hand
side of this equation will be represented as the sum of a series
with unknown coefficients. In Sec.~\ref{Sec DE} we find these
coefficients, solve the differential equation, and get
\eqref{main}.

\section{``Partial densities'' and continuity of $G(t)$} \label{Sec Diff}
We will use the {\bf bold type} for multi-dimensional variables.
Let us indicate the dimension of these variables with subscript,
e.g., $\mathbf{x}_n \in \R^n$; we will omit these subscripts as
often as possible. The coordinates of $\mathbf{x}_n$ will be
always denoted by $x_1, \dots, x_n$. By $\mathbf{0} =
\mathbf{0}_n$ and $\mathbf{1} = \mathbf{1}_n$ denote $(0, \dots
,0)^\top \in \R^n$ and $(1, \dots ,1)^\top \in \R^n$,
respectively.

Let $X_i$ be increments of the exponential random walk $S_i = X_1
+ \dots + X_i$. By definition, $X_i$ are i.i.d. standard
exponential random variables. Put $Y_k:=\frac{2}{k(k+1)}
\sum_{i=1}^k S_i$; clearly, we have $Y_k=\frac{2}{k(k+1)}
\sum_{i=1}^k (k-i+1) X_i$. Hence $\mathbf{Y}_n=(Y_1, \dots,
Y_n)^\top$ is a linear function of $\mathbf{X}_n=(X_1, \dots,
X_n)^\top$, that is $\mathbf{Y}_n =A_n \mathbf{X}_n$, where
$$A_n:=
\begin{pmatrix}
1 & 0 & 0 & \ldots & 0\\
2/3 & 1/3 & 0& \ldots & 0\\
3/6 & 2/6 & 1/6& \ldots & 0\\
\vdots & \vdots & \vdots & \ddots & 0\\
\frac{2n}{n(n+1)} & \frac{2(n-1)}{n(n+1)} & \frac{2(n-2)}{n(n+1)}
& \ldots & \frac{2}{n(n+1)}
\end{pmatrix}.$$
This matrix is nonsingular; by $L_n:= A_n^{-1}$ denote the inverse
matrix. As far as for every $k \ge 3$,
\begin{eqnarray*}
X_k&=&S_k - S_{k-1} = \biggl(\sum_{i=1}^k S_i - \sum_{i=1}^{k-1}
S_i \biggr) - \biggl( \sum_{i=1}^{k-1} S_i - \sum_{i=1}^{k-2} S_i
\biggr) \\
&=& \sum_{i=1}^k S_i  - 2 \sum_{i=1}^{k-1} S_i + \sum_{i=1}^{k-2}
S_i = \frac{k(k+1)}{2} Y_k  - (k-1)k \, Y_{k-1} +
\frac{(k-2)(k-1)}{2} Y_{k-2},
\end{eqnarray*}
we conclude that
$$L_n=
\begin{pmatrix}
1 & 0 & 0 & 0 & \ldots & 0 & 0 & 0\\
-2 & 3 & 0& 0 & \ldots & 0 & 0 & 0\\
1 & -6 & 6& 0&  \ldots & 0 & 0 & 0\\
0 & 3 & -12 & 10 & \ldots & 0 & 0 & 0\\
\vdots & \vdots & \vdots & \vdots & \ddots & \vdots & \vdots & \vdots \\
0 & 0 & 0 & 0 & \ldots & l_{n-2} & 0 & 0\\
0 & 0 & 0 & 0 & \ldots & -2 l_{n-2} & l_{n-1} & 0 \\
0 & 0 & 0 & 0 & \ldots & l_{n-2} & -2 l_{n-1} & l_n
\end{pmatrix}, \qquad l_k:= \frac{k(k+1)}{2}.$$
This matrix has three nonzero diagonals; note that the sum of
elements of each row equals one and the sum of elements of each
column except the last two equals zero.

The distribution of $\mathbf{Y}_n =A \mathbf{X}_n$ is concentrated
on the $n$-dimensional cone $\{\mathbf{y}=\mathbf{y}_n: \, L
\mathbf{y} \ge \mathbf{0} \} \subset \R^n_+$, because $L
\mathbf{Y}_n = \mathbf{X}_n \ge \mathbf{0}$ ($\mathbf{x} \ge
\mathbf{y}$ denotes coordinate-wise inequalities). For the density
$p_{\mathbf{Y}_n}(\mathbf{y}) = p(\mathbf{y})$ of $\mathbf{Y}_n$,
we have $p_{\mathbf{Y}_n}(\mathbf{y}) = |\det A^{-1}|
p_{\mathbf{X}_n}(A^{-1} \mathbf{y}) = \frac{n! (n+1)!}{2^n}
\exp\{- (L\mathbf{y} ,\mathbf{1})\} \,  \I_{ \{ L\mathbf{y} \ge
\mathbf{0} \}}$. For each of $n-2$ first columns of $L$, the sum
of its elements equals zero; the sum of elements in the columns
$n-1$ and $n$ equals $-l_{n-1}=-\frac{(n-1)n}{2}$ and
$l_n=\frac{n(n+1)}{2}$, respectively. Thus we conclude
\begin{equation}
\label{density}
p(\mathbf{y}_n) = \frac{n! (n+1)!}{2^n} e^{- \frac{n(n+1)}{2} y_n
+ \frac{(n-1)n}{2} y_{n-1}} \cdot \I_{ \{ L\mathbf{y} \ge
\mathbf{0} \}}.
\end{equation}
This is not surprise that the density of $\mathbf{Y}_n=(Y_1,
\dots, Y_n)^\top$ depends only on $Y_{n-1}$ and $Y_n$. In fact, we
have $\frac{n(n+1)}{2} Y_n - \frac{(n-1)n}{2} Y_{n-1} = S_n$, and
the density of independent exponential random variables
$\mathbf{X}_n=(X_1, \dots, X_n)^\top$ depends on $S_n$ only.

Now we write
\begin{equation}
\label{G_n = min}
G_n(t)= \P \Bigl \{ \min \limits_{1 \le k \le n } Y_k \ge t \Bigr
\} = \sum_{k=1}^n \P \Bigl \{ \min \limits_{1 \le i \le n} Y_i =
Y_k, Y_k \ge t \Bigr \},
\end{equation}
and for any $1 \le k \le n$ introduce
\begin{equation}
\label{part dens def}
g_n^{(k)}(t) := \int \limits_{ \{ \mathbf{y} = \mathbf{y}_n: \,
\mathbf{y} \ge t\mathbf{1}, \, y_k=t \}} p(\mathbf{y}) d
\lambda_{n-1} (\mathbf{y}),
\end{equation}
where $\lambda_{n-1}$ denotes the Lebesgue measure on the
$(n-1)$-dimensional set $\bigl\{ \mathbf{y}=\mathbf{y}_n: \,
\mathbf{y} \ge t \mathbf{1}, \, y_k=t \bigr\}$. Then by the
Lebesgue theorem, for a.e. $t$,
\begin{equation}
\label{Leb diff}
g_n^{(k)}(t)= - \P \Bigl \{ \min \limits_{1 \le i \le n} Y_i =
Y_k, Y_k \ge t \Bigr \}',
\end{equation}
where the derivative exists a.e. Consequently, $G_n(t)$ is
differentiable a.e., and for almost every $t$, we have
\begin{equation}
\label{G_n'(t) =}
G_n'(t)=- \sum_{k=1}^n g_n^{(k)}(t).
\end{equation}

Let us call $g_n^{(k)}(t)$ the $k$th {\it partial density} of the
random variable $\min \limits_{1 \le k \le n } Y_k$. In the next
section we will show that partial densities are continuous, hence
\eqref{Leb diff} and \eqref{G_n'(t) =} hold for {\it every} $t$
and $G_n(t)$ is differentiable.

Finally, making the change of variables $\mathbf{y} = t
(\mathbf{z} + \mathbf{1})$ in \eqref{part dens def} and using
\eqref{density}, we get
\begin{equation}
\label{part dens}
g_n^{(k)}(t) = \frac{n! (n+1)!}{2^n} t^{n-1} e^{-nt} \int
\limits_{ \{ \mathbf{z} = \mathbf{z}_n: \, L\mathbf{z} \ge
-\mathbf{1}, \, \mathbf{z} \ge \mathbf{0}, \, z_k=0 \}} e^{-
\frac{n(n+1)}{2} tz_n + \frac{(n-1)n}{2} tz_{n-1}} d \lambda_{n-1}
(\mathbf{z}).
\end{equation}
Indeed, the inequality $L\mathbf{y} \ge \mathbf{0}$ transforms to
$L\mathbf{z} \ge \mathbf{-1}$ because for each row of $L$, the sum
of its elements equals one.

We finish the section with the following statement.
\begin{prop}
\label{UNIF G_n}
The functions $G_n(t)$ are continuous. For every $\varepsilon \in
(0, 1)$, $G_n(t)$ converge to $G(t)$ uniformly on $[0, 1 -
\varepsilon]$. The function $G(t)$ is continuous on $[0,1)$.
\end{prop}
\begin{proof}
As $$G_n(t) = \P \biggl \{ \min \limits_{1 \le k \le n}
\sum_{i=1}^k \bigl( S_i - it \bigr) \ge 0 \biggr \},$$ we can
write
\begin{eqnarray*}
G_n(t) - G(t)  &=& \P \biggl \{ \min_{1 \le k \le n } \sum_{i=1}^k
(S_i - i t) \ge 0, \inf_{k > n} \sum_{i=1}^k (S_i -
i t) < 0 \biggr \} \\
&<& \P \biggl \{ \inf_{k > n} \sum_{i=n+1}^k (S_i -
i t) < 0 \biggr \} \\
&<& \P \biggl \{ \exists \, i > n : S_i - i t < 0 \biggr \}.
\end{eqnarray*}
Then $$\sup \limits_{0 \le t \le 1- \varepsilon} \bigl| G_n(t) -
G(t) \bigr| < \sup \limits_{0 \le t \le 1- \varepsilon} \P \biggl
\{ \inf \limits_{i > n} \frac{S_i}{i} < t \biggr \} =\P \biggl \{
\inf \limits_{i > n} \frac{S_i}{i} < 1- \varepsilon \biggr \},$$
and the last expression tends to zero by the strong law of large
numbers.

By \eqref{part dens def} and \eqref{G_n'(t) =}, $G_n(t)$ are
(absolutely) continuous; hence $G(t)$ is continuous on $[0,1)$ as
a uniform limit of continuous functions.
\end{proof}

\section{Properties of ``partial densities'' and differentiability of $G(t)$}
\label{Sec Prop} Here we prove several important properties of
partial densities $g_n^{(k)}(t)$. Let us first state the auxiliary
Lemmata~\ref{PD1},~\ref{CHAINING},~and~\ref{LAST}, and then prove
the differentiability of $G(t)$ in Proposition~\ref{UNIF G_n'}.
The lemmata will be proved afterwards.

\begin{lem}
\label{PD1}
For every $n \ge 1$, $$g_{n+1}^{(1)} (t) = G_n(t) e^{-t}.$$
\end{lem}

\begin{lem}[{\bf Chaining property}]
\label{CHAINING} For every $n \ge 2$ and $1 \le k \le n-1$,
$$g_{n}^{(k)} (t) = c_k \bigl( te^{-t} \bigr)^{k-1}
g_{n-k+1}^{(1)} (t),$$ where $\{ c_k \}_{k \ge 1}$ are some
positive constants. These constants satisfy $c_k = O \bigl(
\sqrt{k} \, e^k \bigr)$.
\end{lem}

\begin{lem}\label{LAST}
The functions $g_n^{(n)} (t)$ are continuous. For every
$\varepsilon \in (0,1)$, $g_n^{(n)} (t)$ converge to $0$ uniformly
on $[0, 1 - \varepsilon]$.
\end{lem}

\begin{prop}\label{UNIF G_n'}
The functions $G_n(t)$ are continuously differentiable. For every
$\varepsilon \in (0, 1)$, $G_n'(t)$ converge to $-G(t) e^{-t}
\sum_{k=1}^\infty c_k \bigl( te^{-t} \bigr)^{k-1}$ uniformly on
$[0, 1 - \varepsilon]$. The function $G(t)$ is continuously
differentiable on $[0,1)$, and
$$G'(t)=-G(t) e^{-t} \sum_{k=1}^\infty c_k \bigl( te^{-t}
\bigr)^{k-1}, \qquad t \in [0,1).$$
\end{prop}

\begin{proof}[{\bf Proof of Proposition~\ref{UNIF G_n'}}]
By \eqref{G_n = min}, \eqref{part dens def}, and \eqref{Leb diff},
we get continuous differentiability of $G_n(t)$ if we show that
partial densities are continuous. The last partial density
$g_n^{(n)}(t)$ is continuous by Lemma~\ref{LAST}. Using
Lemmata~\ref{PD1}~and~\ref{CHAINING},
\begin{equation}
\label{noname1}
g_n^{(k)}(t) = c_k \bigl( te^{-t} \bigr)^{k-1} G_{n-k}(t) e^{-t},
\qquad 1 \le k \le n-1;
\end{equation}
but $G_{n-k}(t)$ are continuous by Proposition~\ref{UNIF G_n},
thus $g_n^{(k)}(t)$ are continuous. Note that now we know that
\eqref{Leb diff} and \eqref{G_n'(t) =} are true for {\it every}
$t$.

Further, $\bigl(t e^{-t} \bigr)' = (1-t)e^{-t} \ge 0$ on $[0,1]$,
thus $t e^{-t} < 1\cdot e^{-1} = e^{-1}$ for $t \in [0, 1)$.
Therefore in view of the estimate on the rate of growth of $c_k$
from Lemma~\ref{CHAINING}, $\sum_{k=1}^\infty c_k \bigl( te^{-t}
\bigr)^{k-1}$ converges on $[0, 1)$. Clearly, this convergence is
uniform on $[0, 1 - \varepsilon]$.

Recall that $1 \ge G_n(t) \searrow G(t) \ge 0$ for every $t$. Then
by \eqref{G_n'(t) =} (which holds for every $t$) and
\eqref{noname1},
\begin{eqnarray*}
&& \biggl| -G(t) e^{-t} \sum_{k=1}^\infty c_k \bigl( te^{-t}
\bigr)^{k-1}-G_n'(t) \biggr| \\
&\le&  e^{-t} \sum_{k=1}^{n-1} c_k \bigl( te^{-t} \bigr)^{k-1}
\bigl( G_{n-k}(t) - G(t) \bigr) + e^{-t} G(t) \sum_{k=n}^\infty
c_k \bigl( te^{-t} \bigr)^{k-1} + g_n^{(n)}(t) \\
&\le&  \bigl( G_{n/2}(t) - G(t) \bigr) \sum_{k=1}^{n/2} c_k \bigl(
te^{-t} \bigr)^{k-1}  + \sum_{k= n/2 + 1 }^\infty c_k \bigl(
te^{-t} \bigr)^{k-1} + g_n^{(n)}(t).
\end{eqnarray*}
The last expression tends to zero uniformly in $t \in [0, 1 -
\varepsilon]$. Indeed, for the third term, use Lemma~\ref{LAST};
the second one is a remainder of the uniformly converging series;
and the first term tends to zero by Proposition~\ref{UNIF G_n}.
\end{proof}

\begin{proof}[{\bf Proof of Lemma~\ref{PD1}}]
The $Y_1=X_1$ is a standard exponential random variable, therefore
\begin{eqnarray*}
g_{n+1}^{(1)} (t) &=& - \P \Bigl \{ \min \limits_{1 \le k \le n+1}
Y_k = Y_1, Y_1 \ge t \Bigr \}' = - \P \Bigl \{ \min \limits_{2 \le
k \le n+1} Y_k \ge Y_1, Y_1 \ge t \Bigr \}'\\
&=& - \left( \int_t^{\infty} \E \Bigl \{ \min \limits_{2 \le k \le
n+1} Y_k \ge Y_1 \Bigl| Y_1 = s \Bigr \} \, d \P \bigl\{ Y_1 < s
\bigr\} \right)' \\
&=& \E \Bigl \{ \min \limits_{2 \le k \le n+1} Y_k \ge Y_1 \Bigl|
Y_1 = t \Bigr \} \, e^{-t}
\end{eqnarray*}
for a.e. $t$. After simple transformations
\begin{eqnarray*}
&& \E \Bigl \{ \min \limits_{2 \le k \le n+1} Y_k \ge Y_1 \Bigl|
Y_1 = t \Bigr \} \\
&=& \E \Bigl \{ \forall \, 2 \le k \le n+1, \, \frac{2}{k(k+1)}
\Bigl( k X_1 + \bigl(k-1 \bigr) X_2 + \dots + X_k \Bigr) \ge X_1
\Bigl| X_1 = t \Bigr \} \\
&=& \E \Bigl \{ \forall \, 2 \le k \le n+1, \, \frac{2}{k(k+1)}
\Bigl(\bigl(k-1 \bigr) X_2 + \dots + X_k \Bigr) \ge \Bigl(1 -
\frac{2}{k+1} \Bigr) X_1 \Bigl| X_1 = t \Bigr \},
\end{eqnarray*}
we use the independence of $X_i$ and find $$\E \Bigl \{ \min
\limits_{2 \le k \le n+1} Y_k \ge Y_1 \Bigl| Y_1 = t \Bigr \} = \P
\Bigl \{ \forall \, 2 \le k \le n+1, \, \frac{2}{(k-1)k}
\Bigl(\bigl(k-1 \bigr) X_2 + \dots + X_k \Bigr) \ge t \Bigr \}.$$
The right-hand side equals $G_n(t)$. It remains to note that the
conditional expectation is continuous because $G_n(t)$ is
continuous (by Proposition~\ref{UNIF G_n}), therefore our argument
is true for {\it every} $t$.
\end{proof}

\begin{proof}[{\bf Proof of Lemma~\ref{CHAINING}}]
The case $k=1$ is trivial, we put $c_1:=1$. Now suppose that $2
\le k \le n-1$.

Let us introduce the following notations. For an $l \times m$
matrix $M$, by $M^{\{i_1, \dots, i_r ; j_1, \dots, j_s \}}$ (where
$r \le l$, $s \le m$ and $1 \le i_1 < \dots < i_r \le l$, $1 \le
j_1 < \dots < j_s \le m$) denote the $(l-r) \times (m-s)$ matrix
obtained from $M$ by deleting the rows $i_1, \dots, i_r$ and the
columns $j_1, \dots, j_s$. For multi-dimensional variables, we
will use the analogous notation.

Consider the integration set $\bigl\{ \mathbf{y} = \mathbf{y}_n:
\, L \mathbf{y} \ge -\mathbf{1}, \, \mathbf{y} \ge \mathbf{0}, \,
y_k=0 \bigr \}$ from \eqref{part dens}. We claim that the first
$k-1$ and the last $n-k$ coordinates of any element of this set
satisfy independent constraints. Precisely, for $2 \le k \le n-1$,
\begin{eqnarray} \label{Fubini}
&& \bigl\{ \mathbf{y}=\mathbf{y}_n: \, L\mathbf{y} \ge
-\mathbf{1}, \, \mathbf{y} \ge \mathbf{0}, \, y_k=0 \bigr \} \notag \\
&=& \bigl\{ \mathbf{y}=\mathbf{y}_{k-1}: \, L_k^{\{\varnothing;
k\}} \mathbf{y} \ge -\mathbf{1}_{k-1}, \, \mathbf{y} \ge
\mathbf{0}_{k-1} \bigr \} \times \bigl\{ 0 \bigr\}\\
&\times& \bigl\{ \mathbf{y}=\mathbf{y}_{n-k}: \, L_n^{\{1, \dots,
k+1 ;1, \dots, k \}}\mathbf{y} \ge -\mathbf{1}_{n-k}, \,
\mathbf{y} \ge \mathbf{0}_{n-k} \bigr \} \notag;
\end{eqnarray}
here $L_k^{\{\varnothing; k\}}$ is $k \times (k-1)$ matrix
consisting of the elements from the top left ``corner'' of $L_n$
and $L_n^{\{1, \dots, k+1 ;1, \dots, k \}}$ is the $(n-k-1) \times
(n-k)$ matrix consisting of the elements from the bottom right
``corner'' of $L_n$.

Take an $\mathbf{y} \in \bigl\{ \mathbf{y}=\mathbf{y}_n: \,
L\mathbf{y} \ge -\mathbf{1}, \, \mathbf{y} \ge \mathbf{0}, \,
y_k=0 \bigr \}$, then $\mathbf{y}^{\{ k\}}$ satisfies
$L_n^{\{\varnothing; k\}} \mathbf{y}^{\{k\}} \ge
-\mathbf{1}_{n-1}$ and $\mathbf{y}^{\{k\}} \ge \mathbf{0}_{n-1}$.
In view of $\mathbf{y}^{\{k\}} \ge \mathbf{0}_{n-1}$, the
$(k+1)$th of the inequalities $L_n^{\{\varnothing; k\}}
\mathbf{y}^{\{k\}} \ge -\mathbf{1}_{n-1}$, namely,
$\frac{(k-1)k}{2} y_{k-1} + \frac{(k+1)(k+2)}{2} y_{k+1} \ge -1$,
holds automatically. Therefore we can delete the row $k+1$ from $
L_n^{\{\varnothing; k\}}$, that is
$$ \bigl \{ L_n^{\{\varnothing; k\}} \mathbf{y}^{\{k\}} \ge
-\mathbf{1}_{n-1}, \, \mathbf{y}^{\{k\}} \ge \mathbf{0}_{n-1}
\bigr \} = \bigl \{ L_n^{\{k+1; k\}} \mathbf{y}^{\{k\}} \ge
-\mathbf{1}_{n-1}, \, \mathbf{y}^{\{k\}} \ge \mathbf{0}_{n-1}
\bigr \}.$$ Since $L_n$ has three nonzero diagonals,
$$L_n^{\{k+1; k\}} =
\begin{pmatrix}
L_k^{\{\varnothing; k\}} & 0 \\
0 & L_n^{\{1, \dots, k+1 ;1, \dots, k \}} \\
\end{pmatrix}$$
is a block matrix. This implies that the constraints for $y_1,
\dots, y_{k-1}$ and $y_{k+1}, \dots, y_n$ are independent, i.e.,
\eqref{Fubini} holds true.

Now define
\begin{eqnarray}
\label{notat P v c}
&P_m := \bigl\{ \mathbf{y}=\mathbf{y}_{m-1}: \,
L_m^{\{\varnothing; m\}}\mathbf{y} \ge -\mathbf{1}_{m-1}, \,
\mathbf{y} \ge \mathbf{0}_{m-1} \bigr \},& \quad m \ge 2 \notag \\
&v_m := \lambda_{m-1} (P_m),&\\
&c_m := 2^{-m} m!(m+1)! v_m,& \notag
\end{eqnarray}
and combine \eqref{Fubini} with Fubini's theorem to rewrite
\eqref{part dens} in a simpler form. For $k=n-1$ we get
\begin{equation}
\label{part d = n-1}
g_n^{(n-1)}(t) = \frac{n! (n+1)!}{2^n} t^{n-1} e^{-nt} v_{n-1}
\int \limits_{ \{ y \ge 0 \}} e^{- \frac{n(n+1)}{2} t y} d
\lambda_1 (y),
\end{equation}
and for $2 \le k \le n-2$ we get
\begin{eqnarray}
\label{part d < n-1}
g_n^{(k)}(t) &=& \frac{n! (n+1)!}{2^n} t^{n-1} e^{-nt} v_k \\
& \times& \int \limits_{ \bigl\{ \mathbf{y}=\mathbf{y}_{n-k}: \,
L_n^{\{1, \dots, k+1 ;1, \dots, k \}}\mathbf{y} \ge
-\mathbf{1}_{n-k}, \, \mathbf{y} \ge \mathbf{0}_{n-k} \bigr \} }
e^{- \frac{n(n+1)}{2} ty_{n-k} + \frac{(n-1)n}{2} ty_{n-k-1}} d
\lambda_{n-k} (\mathbf{y}) \notag.
\end{eqnarray}

For the simpler case $k=n-1$, we integrate in \eqref{part d = n-1}
and find
\begin{equation}
\label{next to last}
g_n^{(n-1)}(t) = c_{n-1} t^{n-2} e^{-nt}.
\end{equation}
But by Lemma~\ref{PD1}, $g_2^{(1)} (t) = e^{-2t}$, and there is
nothing to prove.

For the harder case $2 \le k \le n-2$, we use \eqref{part dens}
and write
\begin{eqnarray}
\label{long}
g_{n-k+1}^{(1)}(t) &=& \frac{(n-k+1)! (n-k+2)!}{2^{n-k+1}} t^{n-k}
e^{-(n-k+1)t} \\
& \times & \int \limits_{ \{ \mathbf{z} = \mathbf{z}_{n-k+1}: \,
L\mathbf{z} \ge -\mathbf{1}, \, \mathbf{z} \ge \mathbf{0}, \,
z_1=0 \}} e^{- \frac{(n-k+1)(n-k+2)}{2} tz_{n-k+1} +
\frac{(n-k)(n-k+1)}{2} tz_{n-k}} d \lambda_{n-k} (\mathbf{z})
\notag.
\end{eqnarray}
Take an element $\mathbf{z}$ of the integration set, then
$\mathbf{z}^{\{ 1\}}$ satisfies $L_{n-k+1}^{\{\varnothing; 1\}}
\mathbf{z}^{\{1\}} \ge -\mathbf{1}_{n-k}$ and $\mathbf{z}^{\{1\}}
\ge \mathbf{0}_{n-k}$. But the first of the inequalities
$L_{n-k+1}^{\{\varnothing; 1\}} \mathbf{z}^{\{1\}} \ge
-\mathbf{1}_{n-k}$, that is $0 \ge -1$, is always true, while the
second one, $3 z_2 \ge -1$, follows from $\mathbf{z}^{\{1\}} \ge
\mathbf{0}_{n-k}$. Therefore
\begin{eqnarray*}
&&\bigl \{ \mathbf{z} = \mathbf{z}_{n-k+1}: \, L_{n-k+1}\mathbf{z}
\ge -\mathbf{1}, \, \mathbf{z} \ge \mathbf{0}, \, z_1=0 \bigr \} \\
&=& \bigl \{ 0 \bigr \} \times \bigl \{ \mathbf{z} =
\mathbf{z}_{n-k}: \, L_{n-k+1}^{ \{1,2 ; 1\}} \mathbf{z} \ge
-\mathbf{1}_{n-k}, \, \mathbf{z} \ge \mathbf{0}_{n-k} \bigr \},
\end{eqnarray*}
and transforming the integral in \eqref{long}, we get
\begin{eqnarray}
\label{integral}
g_{n-k+1}^{(1)}(t) &=& \frac{(n-k+1)! (n-k+2)!}{2^{n-k+1}} t^{n-k}
e^{-(n-k+1)t} \\
& \times & \int \limits_{ \{ \mathbf{z} = \mathbf{z}_{n-k}: \,
L_{n-k+1}^{ \{1,2 ; 1\}} \mathbf{z} \ge -\mathbf{1}_{n-k}, \,
\mathbf{z} \ge \mathbf{0}_{n-k} \}} e^{- \frac{(n-k+1)(n-k+2)}{2}
tz_{n-k} + \frac{(n-k)(n-k+1)}{2} tz_{n-k-1}} d \lambda_{n-k}
(\mathbf{z}) \notag.
\end{eqnarray}

Compare \eqref{part d < n-1} and \eqref{integral}. The integrals
in these formulas do not look nice, but the point is that one
could be obtained from the other by a very simple change of
variables. Recall that
$$L_{n-k+1}^{ \{1,2 ; 1\}} =
\begin{pmatrix}
-6 & 6& 0&  \ldots & 0 & 0\\
 3 & -12 & 10 & \ldots & 0 & 0\\
\vdots & \vdots & \vdots & \ddots & \vdots & \vdots\\
0 & 0 & 0 & \ldots & l_{n-k} & 0 \\
0 & 0 & 0 & \ldots &  -2l_{n-k} & l_{n-k+1}
\end{pmatrix}$$
and
$$L_n^{\{1, \dots, k+1 ;1, \dots, k \}} =
\begin{pmatrix}
-2l_{k+1} & l_{k+2} & 0& \ldots & 0 & 0\\
l_{k+1} & -2l_{k+2} & l_{k+3}& \ldots & 0 & 0\\
\vdots & \vdots & \vdots & \ddots & \vdots & \vdots \\
0 & 0 & 0 & \ldots &  l_{n-1} & 0 \\
0 & 0 & 0 & \ldots & -2l_{n-1} & l_n
\end{pmatrix}.$$
By $J_m^{s \mapsto i}$, where $i, s, m \ge 1$, we denote the $m
\times m$ diagonal matrix with the elements $\frac{l_i}{l_s},
\frac{l_{i+1}} {l_{s+1}}, \dots, \frac{l_{i+m-1}}{l_{s+m-1}}$ on
the diagonal (counting from the top left corner). Then
\begin{equation}
\label{L = LJ}
L_{n-k+1}^{ \{1,2 ; 1\}} = L_n^{\{1, \dots, k+1 ;1, \dots, k \}}
J_{n-k}^{k+1 \mapsto 2}
\end{equation}
because the right-sided multiplication by $J_{n-k}^{k+1 \mapsto
2}$ multiplies the first column of $L_n^{\{1, \dots, k+1 ;1,
\dots, k \}}$ by $\frac{l_2}{l_{k+1}} = \frac{3}{l_{k+1}}$, the
second column by $\frac{l_3}{l_{k+2}} = \frac{6}{l_{k+2}}$, etc.

Hence the change of variables $\mathbf{z} = \bigl( J_{n-k}^{k+1
\mapsto 2} \bigr)^{-1} \mathbf{y}$ transforms the integral from
\eqref{integral} to the integral from \eqref{part d < n-1} times
$\bigl| \det \bigl( J_{n-k}^{k+1 \mapsto 2} \bigr)^{-1} \bigr|$.
But
\begin{equation}
\label{det J}
\det J_m^{s \mapsto i} = \frac{(s-1)!s! (i+m-1)!(i+m)!}{(i-1)!i!
(s+m-1)! (s+m)!},
\end{equation}
 and from \eqref{part d < n-1} and \eqref{integral},
$$\frac{g_n^{(k)}(t)}{g_{n-k+1}^{(1)}(t)} = \frac{\frac{n! (n+1)!}{2^n}
t^{n-1} e^{-nt} v_k}{ \det \bigl( J_{n-k}^{k+1 \mapsto 2}
\bigr)^{-1} \cdot \frac{(n-k+1)! (n-k+2)!}{2^{n-k+1}} t^{n-k}
e^{-(n-k+1)t}} = \frac{k!(k+1)! v_k}{2^k} t^{k-1}e^{-(k-1)t}.$$

It remains to check that $c_k = O \bigl( \sqrt{k} \, e^k \bigr)$
to finish the proof of Lemma. By \eqref{next to last},
$g_{k+1}^{(k)}(t) = c_k t^{k-1} e^{-(k+1)t}$, and integrating from
$0$ to $\infty$, we get $\P \Bigl \{ \min \limits_{1 \le i \le
k+1} Y_i = Y_k \Bigr \} = c_k \frac{(k-1)!}{(k+1)^k}$. Then $c_k <
\frac{(k+1)^k}{(k-1)!}$, and by Stirling's formula, $$c_k \le C
\frac{(k+1)^k e^{k-1}}{(k-1)^{k-1} \sqrt{2 \pi (k-1)}} =
C\left(\frac{k+1}{k-1}\right)^{k-1} \frac{k+1 }{\sqrt{2 \pi
(k-1)}} e^{k-1} = O \bigl( \sqrt{k} \, e^k \bigr).$$
\end{proof}

\begin{proof}[{\bf Proof of Lemma~\ref{LAST}}]
Using notations \eqref{notat P v c}, formula \eqref{part dens}
could be written as
\begin{equation}
\label{last part dens}
g_n^{(n)}(t) = \frac{n! (n+1)!}{2^n} t^{n-1} e^{-nt} \int
\limits_{ P_n} e^{\frac{(n-1)n}{2} ty_{n-1}} d \lambda_{n-1}
(\mathbf{y}_{n-1}).
\end{equation}
Recall that, first, $\mathbf{y}_{n-1} \in P_n$ implies
$\mathbf{y}_{n-1} \ge \mathbf{0}$, and second, from the proof of
Lemma~\ref{CHAINING} we know that $P_n \subset \R^{n-1}$ has
finite volume. Then $P_n$ is bounded since it is an intersection
of half-spaces; note that we give another proof of boundedness of
$P_n$ while proving Lemma~\ref{V} from the next section.

Now it is clear that $g_n^{(n)}(t)$ is continuous. Further, as far
as $y_{n-1} \ge 0$ for $\mathbf{y}_{n-1} \in P_n$, the integral
from \eqref{last part dens} increases in $t$. By $\bigl( t^{n-1}
e^{-nt} \bigr )' = t^{n-2} e^{-nt} \bigl( (n-1) - nt ) \bigr)$, we
conclude that $g_n^{(n)}(t)$ increases (in $t$) at least on $[0,
\frac{n-1}{n}]$. Now suppose $n$ is such that $\frac{n-1}{n} \ge 1
- \varepsilon/2$; then $g_n^{(n)}$ is increasing on $[0,
1-\varepsilon/2]$, and in view of \eqref{Leb diff} (which holds
for every $t$), $$\sup \limits_{0 \le t \le 1- \varepsilon}
g_n^{(n)} (t) = g_n^{(n)} (1 - \varepsilon) \le
\frac{2}{\varepsilon} \int \limits_{1 - \varepsilon}^{1 -
\varepsilon/2} g_n^{(n)} (s) ds \le \frac{2}{\varepsilon} \, \P
\Bigl \{ Y_n \le 1 - \varepsilon/2 \Bigr \}.$$ The last expression
tends to zero as $n \to \infty$. Indeed, $Y_n = \frac{2}{n(n+1)}
\sum_{i=1}^n S_i \to 1$ a.s. because $S_n \to 1$ a.s. by the
strong law of large numbers.
\end{proof}

\section{The differential equation for $G(t)$} \label{Sec DE}
By Proposition~\ref{UNIF G_n'}, we know that $G(t)$ satisfies
\begin{equation}
\label{G'(t) =}
\begin{cases}
G'(t)=-G(t) t^{-1} f\bigl( te^{-t} \bigr), \quad t \in [0,1)\\
G(0)=1,
\end{cases}
\end{equation}
where $f(x):= \sum_{k=1}^\infty c_k x^k$ is the generating
function of $c_k$. By Lemma~\ref{CHAINING}, this series converges
for $|x| < e^{-1}$. Let us formulate Lemma~\ref{V}, which is
indispensable for finding $f(x)$, and then find $G(t)$ in
Proposition~\ref{FINAL}. The lemma will be proved afterwards.
\begin{lem}
\label{V}
For $n \ge 2$, we have
$$c_n = \frac{n(n+1)}{n-1} \sum_{k=1}^{n-1} \frac{c_k
c_{n-k}}{(k+1)(n-k+1)},$$ and $c_1=1$.
\end{lem}

\begin{prop} \label{FINAL}
The function $G(t)$ satisfies the differential equation
\begin{equation}
\label{diff eq G}
\begin{cases}
G'(t)= \frac{t-2}{2(1-t)} G(t), \quad t \in [0,1)\\
G(0)=1,
\end{cases}
\end{equation}
which has a unique on $[0,1)$ solution $\sqrt{1-t} e^{-t/2}$.
\end{prop}
\begin{proof}[{\bf Proof of Proposition~\ref{FINAL}}]
Define the variables $b_n:= \frac{c_n}{n+1}$. The generating
function $h(x):= \sum_{k=1}^\infty b_k x^k$ of these variables
satisfies
\begin{equation}
\label{gen func}
\bigl ( x h(x) \bigr) ' = f(x), \qquad |x| < e^{-1}
\end{equation}
and $h(0) = 0$; we recall that the sum of a power series could be
differentiated termwise inside the circle of its convergence.

Using Lemma~\ref{V}, we find that $b_n = \frac{n}{n-1}
\sum_{k=1}^{n-1} b_k b_{n-k}$ for $n \ge 2$ and  $b_1=1/2$. Then
for $|x| < e^{-1}$, $$h^2(x) = \bigl( b_1 x + b_2 x^2 + \dots
\bigr)^2 = \sum_{n = 2}^\infty \sum_{k=1}^{n-1} b_k b_{n-k} x^n =
\sum_{n = 1}^\infty \frac{n-1}{n} b_n x^n = \sum_{n = 1}^\infty
b_n x^n - \sum_{n =1}^\infty \frac{b_n}{n} x^n,$$ and by
differentiation,
$$2h(x)h'(x) = h'(x) - \frac{h(x)}{x}, \qquad |x| < e^{-1}.$$

Now we have $$\begin{cases}
\frac{2h-1}{h}h' = - \frac{1}{x}, \qquad |x| < e^{-1}\\
h'(0)=b_1=1/2,
\end{cases}$$
and by taking into account that $h > 0$ for $x >0$,
$$2h - \ln h = -\ln x + C, \qquad 0 \le x < e^{-1};$$ then
$$\ln 2h - 2h = \ln x + C', \qquad 0 \le x < e^{-1},$$
$$2h e^{-2h} = C'' x, \qquad 0 \le x < e^{-1}.$$ Using $h'(0)=1/2$,
we find that $C''=1$. The function $q(x):= x e^{-x}$ is invertible
on $[0, 1]$, and we obtain $$h(x) = \frac{q^{-1}(x)}{2}, \qquad 0
\le x < e^{-1}.$$

Finally, by \eqref{gen func},
\begin{eqnarray*}
f(te^{-t}) &=& f(q(t))=h(q(t)) + q(t) h'(q(t)) = h(q(t)) +
q(t)\frac{h(q(t))'}{q'(t)}\\
&=& \frac{t}{2} + t e^{-t} \frac{1/2}{(1-t) e^{-t}} =
\frac{2-t}{2(1-t)}t, \qquad t \in [0,1),
\end{eqnarray*}
and applying \eqref{G'(t) =}, we see that $G(t)$ satisfies
\eqref{diff eq G}.

It remains to solve \eqref{diff eq G}. As far as
$\frac{t-2}{2(1-t)} = - \frac{1}{2} - \frac{1}{2(1-t)}$, we get
$$
\begin{cases}
\ln G(t)= -\frac{t}{2} + \frac{1}{2} \ln (1-t) + C, \quad t \in [0,1)\\
G(0)=1.
\end{cases}
$$
Then $C=0$, and $G(t)= \sqrt{1-t} e^{-t/2}$ on $[0,1)$. Note that
this equality holds on $[0,1]$.
\end{proof}

\begin{proof}[{\bf Proof of Lemma~\ref{V}}]
We recall that $c_n = 2^{-n} n!(n+1)! v_n$ for $n \ge 2$, where
$v_n$ is the volume of $$P_n = \bigl\{
\mathbf{y}=\mathbf{y}_{n-1}: \, L_n^{\{\varnothing; n\}}\mathbf{y}
\ge -\mathbf{1}_{n-1}, \, \mathbf{y} \ge \mathbf{0}_{n-1} \bigr
\},$$ see \eqref{notat P v c}. That is why our goal is to find
$v_n$. As far as $c_1=1$, we define $v_1:=1$ to satisfy $c_1 =
2^{-1} 1!(1+1)! v_1$.

Let us study properties of $P_n$. We will temporary forget about
geometric intuition and use algebraic arguments only. Evidently,
$$P_n = \bigl\{ \mathbf{y}=\mathbf{y}_{n-1}: \, L_n^{\{1; n\}}
\mathbf{y} \ge -\mathbf{1}_{n-1}, \, \mathbf{y} \ge
\mathbf{0}_{n-1} \bigr \}$$ because the first of the inequalities
$\, L_n^{\{\varnothing; n\}} \mathbf{y} \ge -\mathbf{1}_{n-1}$,
namely, $y_1 \ge -1$, follows from $\mathbf{y} \ge
\mathbf{0}_{n-1}$.

We claim that (for every $n \ge 2$) the matrix $L_n^{\{1; n\}}$ is
nonsingular. By $I^\vartriangle_k$ denote the $k \times k$ upper
triangular matrix with all its $\frac{k(k+1)}{2}$ nonzero elements
equal $1$, and denote $I^\triangledown_k :=
(I^\vartriangle_k)^\top$, which is lower triangular. For the
matrix $L_n^{\{1; n\}}$, the sum of elements of each column except
the first and the last ones equals zero. For the matrix
$I^\triangledown_{n-1} L_n^{\{1; n\}}$, the sum of elements of
each column except the first one equals zero. That is why we
easily get
\begin{equation}
\label{2 tri} I^\vartriangle_{n-1} I^\triangledown_{n-1} L_n^{\{1;
n\}} =
\begin{pmatrix}
-n  & 0 & 0 & \ldots & 0 & 0\\
-n+2  &  -3 &0 & \ldots & 0 & 0\\
-n+3  &  0 &-6 & \ldots & 0 & 0\\
\vdots & \vdots & \vdots & \ddots & \vdots & \vdots\\
-2 & 0 & 0 & \ldots & -l_{n-2} & 0 \\
-1 & 0 & 0 & \ldots &  0& -l_{n-1}
\end{pmatrix},
\end{equation}
and thus $\det L_n^{\{1; n\}} = \det (I^\vartriangle_{n-1}
I^\triangledown_{n-1} L_n^{\{1; n\}}) \neq 0$.

Let us now show that $L_n^{\{1; n\}} \mathbf{y} \ge \mathbf{0}$
implies $\mathbf{y} \le \mathbf{0}$. The proof is by induction.
For $n=2$, the statement is trivial. Assume that the statement is
true for an $n \ge 2$; then check that it holds for $n+1$. Since
$I^\vartriangle_n I^\triangledown_n$ has positive elements only,
$L_{n+1}^{\{1; n+1\}} \mathbf{y} \ge \mathbf{0}$ implies
$I^\vartriangle_n I^\triangledown_n L_{n+1}^{\{1; n+1\}}
\mathbf{y} \ge \mathbf{0}$, and by \eqref{2 tri} we conclude that
$y_1 \le 0$. Further, $L_{n+1}^{\{1; n+1\}} \mathbf{y} \ge
\mathbf{0}$ implies $L_{n+1}^{\{1,2; 1,n+1\}} \mathbf{y}^{\{1\}}
\ge \mathbf{0}_{n-1}$. In fact, we just get rid of the first of
the inequalities $L_{n+1}^{\{1; n+1\}} \mathbf{y} \ge \mathbf{0}$
and replace $y_1 - 6 y_2 + 6 y_3 \ge 0$, which is the second
inequality, by the less restrictive $ - 6 y_2 + 6 y_3 \ge 0$
(recall that $y_1 \le 0$ !). Then $L_{n+1}^{\{1,2; 1,n+1\}} =
L_n^{\{1; n\}} J_{n-1}^{1 \mapsto 2}$ (see the comment to
analogous statement \eqref{L = LJ} and the definition of $J_m^{s
\mapsto i}$), therefore $L_{n+1}^{\{1,2; 1,n+1\}}
\mathbf{y}^{\{1\}} \ge \mathbf{0}_{n-1}$ is equivalent to
$L_n^{\{1; n\}} \mathbf{z} \ge \mathbf{0}_{n-1}$, where
$\mathbf{z} := J_{n-1}^{1 \mapsto 2} \mathbf{y}^{\{1\}}$. By the
assumption, $\mathbf{z} \le \mathbf{0}_{n-1}$; hence
$\mathbf{y}^{\{1\}} \le \mathbf{0}_{n-1}$ because $\mathbf{z}$ is
obtained from $\mathbf{y}^{\{1\}}$ by the tension with positive
coefficients. Finally, we get $\mathbf{y} \le \mathbf{0}_n$.

By \eqref{2 tri} we easily find that a unique solution of
$L_n^{\{1; n\}}\mathbf{y} = -\mathbf{1}$ is
\begin{equation}
\label{y^*}
\mathbf{y}^*_{n-1}:=\Bigl( \frac{n-1}{2}, \frac{n-2}{3}, \dots,
\frac{2}{n-1}, \frac{1}{n}\Bigr)^\top.
\end{equation}
We see that $\mathbf{y}^* \in P_n$, and since $L_n^{\{1; n\}}
(\mathbf{y} - \mathbf{y}^*) \ge \mathbf{0}_{n-1}$ for every
$\mathbf{y} \in P_n$, we have $\mathbf{y} \le \mathbf{y}^*$.
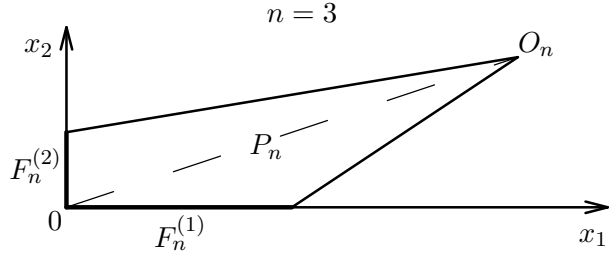
\begin{center}
\begin{texdraw}
\drawdim cm \setunitscale 6

\linewd 0.005 \arrowheadsize l:0.05 w:0.025 \arrowheadtype t:V
\move (0.1 0.1) \avec(1.3 0.1) \move (1.3 0.05) \textref h:R v:T
\htext{$x_1$}

\move (0.1 0.1) \avec(0.1 0.5) \move (0.07 0.47) \textref h:R v:T
\htext{$x_2$}

\move (0.6 0.1) \lvec (1.1 0.433) \lvec (0.1 0.267)

\linewd 0.01 \move (0.1 0.1) \lvec (0.6 0.1) \move (0.1 0.1) \lvec
(0.1 0.267)

\move (0.7 0.55) \textref h:R v:T \htext{$n=3$}

\move (0.09 0.09) \textref h:R v:T \htext{$0$} \move (0.09 0.23)
\textref h:R v:T \htext{$F_n^{(2)}$}

\move (0.41 0.08 ) \textref h:R v:T \htext{$F_n^{(1)}$}

\linewd 0.003 \move (0.1 0.1) \lpatt (0.11 0.14) \lvec (1.1 0.433)
\textref h:L v:B \htext{$O_n$}

\move (0.58 0.26 ) \textref h:R v:T \htext{$P_n$}

\move (0.1 0.6) \move (0.1 0.01)
\end{texdraw}
\end{center}

Now it is clear that $P_n$ is an $(n-1)$-dimensional convex
polyhedron with $2(n-1)$ faces. Denote by $O_n$ the point of
intersection of $n-1$ hyperplanes (faces) $L_n^{\{1;
n\}}\mathbf{y} = -\mathbf{1}_{n-1}$, and put $F_n^{(k)}:= P_n \cap
\bigl\{ \mathbf{y} = \mathbf{y}_{n-1}: y_k=0 \bigr \}$, where $1
\le k \le n-1$, see the figure. Naturally, $O_n =
\mathbf{y}^*_{n-1}$, and by $\mathbf{y}_{n-1} \in P_n
\Longrightarrow \mathbf{y}_{n-1} \le \mathbf{y}^*_{n-1}$, the
$P_n$ is a disjoint union of $n-1$ simpleces with the common
vertex $O_n$ and the bases $F_n^{(k)}$ (to be pedantic, the
simpleces themselves are not disjoint, but their interiors are).
Recalling \eqref{y^*},
\begin{equation}
\label{v_n =}
v_n = \frac{1}{n-1} \sum_{k=1}^{n-1} \frac{n-k}{k+1} \cdot
\lambda_{n-2}(F_n^{(k)}).
\end{equation}
Thus we reduced the problem of finding the volume of $P_n$ to
finding the volumes of $F_n^{(k)}$.

For $2 \le k \le n-2$, we have
\begin{eqnarray*}
F_n^{(k)} &=& \bigl\{ \mathbf{y}=\mathbf{y}_{n-1}: \, L_n^{\{1;
n\}}\mathbf{y} \ge -\mathbf{1}_{n-1}, \, \mathbf{y} \ge
\mathbf{0}_{n-1}, y_k = 0 \bigr \} \\
&=& \bigl\{ \mathbf{y}=\mathbf{y}_{k-1}: \, L_k^{\{1; k\}}
\mathbf{y} \ge -\mathbf{1}_{k-1}, \, \mathbf{y} \ge
\mathbf{0}_{k-1} \bigr \} \times \bigl\{ 0 \bigr \} \\
&\times& \bigl\{ \mathbf{y}=\mathbf{y}_{n-k-1}: \, L_n^{\{1,
\dots, k+1 ;1, \dots, k,n \}} \mathbf{y} \ge -\mathbf{1}_{n-k-1},
\, \mathbf{y} \ge \mathbf{0}_{n-k-1} \bigr \},
\end{eqnarray*}
which is very similar to \eqref{Fubini}. The only difference is
that this statement deals with $L_n^{\{1; n \}}$ while
\eqref{Fubini} is a statement about $L_n$. Therefore the proof is
a verbatim copy of the proof of \eqref{Fubini}. Then
$$F_n^{(k)} = P_k \times \bigl\{ 0 \bigr \} \times \bigl\{
\mathbf{y}=\mathbf{y}_{n-k-1}: \, L_n^{\{1, \dots, k+1 ;1, \dots,
k,n \}} \mathbf{y} \ge -\mathbf{1}_{n-k-1}, \, \mathbf{y} \ge
\mathbf{0}_{n-k-1} \bigr \},$$ and in view of $L_{n-k}^{\{1 ; n-k
\}} = L_n^{\{1, \dots, k+1 ;1, \dots, k,n \}} J_{n-k-1}^{k+1
\mapsto 1}$ (see the comment to analogous statement \eqref{L = LJ}
and the definition of $J_m^{s \mapsto i}$), we make the change of
variables $\mathbf{y} = J_{n-k-1}^{k+1 \mapsto 1}\mathbf{z}$ and
obtain
$$F_n^{(k)} = P_k \times \bigl\{ 0 \bigr \}
\times J_{n-k-1}^{k+1 \mapsto 1} \bigl\{
\mathbf{z}=\mathbf{z}_{n-k-1}: \, L_{n-k}^{\{1 ; n-k \}}
\mathbf{z} \ge -\mathbf{1}_{n-k-1}, \, \mathbf{z} \ge
\mathbf{0}_{n-k-1} \bigr \}.$$ Finally, $F_n^{(k)} = P_k \times
\bigl\{ 0 \bigr \} \times J_{n-k-1}^{k+1 \mapsto 1} P_{n-k}$, and
by \eqref{det J},
\begin{equation}
\label{vol F}
\lambda_{n-2}(F_n^{(k)}) = \frac{k! (k+1)! (n-k-1)! (n-k)!}{(n-1)!
n!} v_k v_{n-k}.
\end{equation}

For $k=1$, we could repeat word by word the argument we used in
the case $2 \le k \le n-2$. At the final step, we get $F_n^{(1)} =
\bigl\{ 0 \bigr \} \times J_{n-2}^{2 \mapsto 1} P_{n-1}$, and
since $v_1=1$, \eqref{vol F} also holds for $k=1$.

For $k = n-1$, it is readily seen that $F_n^{(n-1)} = P_{n-1}
\times \{ 0 \}$, hence \eqref{vol F} holds for $k=n-1$.

Thus \eqref{vol F} is true for $1 \le k \le n-1$. It remains to
take a look at \eqref{v_n =} and \eqref{vol F} to get
$$v_n = \frac{1}{(n-1) (n-1)! n!} \sum_{k=1}^{n-1} \bigl( (n-k)!
k! \bigr)^2 v_k v_{n-k}.$$ The application of
$v_n=\frac{2^n}{n!(n+1)!} c_n$ finishes the proof of Lemma.
\end{proof}

Let us finish the paper proving the equality $$ G(t)= \lim
\limits_{n \to \infty} \P \biggl \{ \min \limits_{1 \le k \le n}
\frac{2n}{k(k+1)} \sum_{i=1}^k U_{i, n}  \ge t \biggr \},$$ where
$U_{i, n}$ are the order statistics of $n$ i.i.d. random variables
uniformly distributed on $[0,1]$.

It is well-known (see, for example, \cite[Sec. 9.1]{Karlin}) that
$$\bigl( U_{1, n}, U_{2, n}, \dots, U_{n, n} \bigr) \stackrel{\D}{=}
\Bigl( \frac{S_1}{S_{n+1}}, \frac{S_2}{S_{n+1}}, \dots,
\frac{S_n}{S_{n+1}} \Bigr)$$ for every $n \ge 1$. Then $$\P \biggl
\{ \min \limits_{1 \le k \le n } \frac{2n}{k(k+1)} \sum_{i=1}^k
U_{i, n} \ge t \biggr \} = \P \biggl \{ \min \limits_{1 \le k \le
n } \frac{2}{k(k+1)} \sum_{i=1}^k S_i \ge \frac{S_{n+1}}{n} \, t
\biggr \},$$ and combining the law of large numbers with the
pointwise convergence of $G_n$ to $G$, for every $\varepsilon \in
(0,1)$, we have
\begin{eqnarray*}
G \bigl((1+\varepsilon) t\bigr) &\le& \varliminf \limits_{n \to
\infty} \P \biggl \{ \min \limits_{1 \le k \le n }
\frac{2n}{k(k+1)} \sum_{i=1}^k U_{i, n} \ge t \biggr \} \\
&\le& \varlimsup \limits_{n \to \infty} \P \biggl \{ \min
\limits_{1 \le k \le n } \frac{2n}{k(k+1)} \sum_{i=1}^k U_{i, n}
\ge t \biggr \} \le G \bigl((1-\varepsilon) t \bigr).
\end{eqnarray*}
We complete the proof proceeding to the limit $\varepsilon
\searrow 0$ and using the continuity of $G(t)$.

\section*{Acknowledgements}
The author is grateful to M. Lifshits for his interest to this
paper.

\bigskip
\begin{tabular}{>{\footnotesize} l}
Vladislav V. Vysotsky \\
Department of Probability Theory and Mathematical Statistics\\
Faculty of Mathematics and Mechanics\\
St. Petersburg State University\\
Bibliotechnaya pl., 2 \\
Stary Peterhof, 198504\\
Russia\\
E-mail: vlad.vysotsky@gmail.com, vysotsky@vv9034.spb.edu\\
\end{tabular}

\end{document}